\documentclass[12pt]{article}
\usepackage{amssymb,amsmath}
\newtheorem{theorem}{Theorem}[section]
\newtheorem{definition}{Definition}[section]

\newtheorem{lemma}{Lemma}[section]
\newtheorem{proposition}{Proposition}[section]
\newtheorem{example}{Example}[section]

\newtheorem{remark}{Remark}[section]
\newtheorem{corollary}{Corollary}[section]

\begin{document}
\begin{center}
{\bf Best Proximity Point Theorems for Asymptotically Relatively Nonexpansive Mappings}\\

\vspace{0.5cm}

{ S. Rajesh{\footnote{This is an Accepted Manuscript of an article published by Taylor and Francis in [Numerical functional analysis and optimization] on [29 Dec 2015], available online: http://wwww.tandfonline.com/[DOI: 10.1080/01630563.2015.1079533].}}  and P. Veeramani }

{Department of Mathematics, Indian Institute of Technology Madras}

{Chennai - 600036, India}

{e-mail: srajeshiitmdt@gmail.com, pvmani@iitm.ac.in}

\end{center}

\begin{abstract}
Let $(A, B)$ be a nonempty bounded closed convex
proximal parallel pair in a nearly uniformly convex Banach
space and $T: A\cup B \rightarrow A\cup B$ be a continuous and
asymptotically relatively
nonexpansive map. We prove that there exists $x \in A\cup B$
such that $\|x - Tx\| = \emph{dist}(A, B)$ whenever
$T(A) \subseteq B$, $T(B) \subseteq A$. Also, we establish that
if $T(A) \subseteq A$ and $T(B) \subseteq B$, then there exist
$x \in A$ and $y\in B$ such that $Tx = x$, $Ty = y$ and
$\|x - y\| = \emph{dist}(A, B)$. We prove the aforesaid results
when the pair $(A, B)$ has the rectangle property and property $UC$.
In case of $A = B$, we obtain, as a particular case of our results, 
the basic fixed point theorem for asymptotically
nonexpansive maps by Goebel and Kirk.
\end{abstract}




\section{Introduction}
\label{intro}
Let $A$ and $B$ be nonempty weakly compact convex subsets
of a Banach space $X$. 
Suppose $T$ is a relatively nonexpansive
self-map on $A\cup B$. Assume that $T$ maps $A$
into $B$ and $B$ into $A$. Then it was proved,
under suitable assumptions, that there exists a point $x$ in $A$
such that the distance between $x$ and $Tx$ is equal to
the distance between $A$ and $B$;
see Theorem 2.1 in \cite{Eldr_2005}.

Also, it was shown, under suitable assumptions,
that there exist $x$ in $A$ and $y$ in $B$ such that
$x$ and $y$ are fixed points of $T$ and
the distance between $x$ and $y$ is equal to the distance between $A$ and $B$
whenever $T$ maps $A$ into $A$
and $B$ into itself; see Theorem 2.2 in \cite{Eldr_2005}.
For recent advancements in the theory of best proximity points, one may refer to
 \cite{Espi_2013,Fern_2014,Sadi_2013,Kosu_2011}
 and the references therein.

If $A$ is equal to $B$, then we obtain the well known Browder-G$\ddot{o}$hde-Kirk fixed point theorem,
as a special case of the above mentioned results.
A good account of metric fixed point theory can be
found in \cite{Bena_Book,GK_Book,Kham_Book}.

 Note that the notion of asymptotically nonexpansive maps,
which is an interesting generalization of nonexpansive maps, was
introduced in \cite{Goeb_1972}.
Also, it was proved in \cite{Goeb_1972} that if $K$ is a nonemtpy bounded
closed convex subset of a uniformly convex Banach space and
$T$ is an asymptotically nonexpansive self map on $K$, then
there exists $x$ in $K$ such that $x$ is a fixed point of $T$.
This result was further extended to nearly uniformly
convex (NUC) Banach spaces by Xu \cite{Xu_1991}.

In this paper, by observing the notions of relatively nonexpansive maps 
and asymptotically nonexpansive maps, we define a concept called
asymptotically relatively nonexpansive maps and study, under suitable
assumptions, the existence of best proximity points of such maps.
Also, we introduce the notion of proximal parallel pairs
having rectangle property and
establish the existence of a
best proximity point for an asymptotically relatively nonexpansive
map on a proximal parallel pair having rectangle property
in a nearly uniformly convex
Banach space. We also prove, under suitable assumptions, a fixed point theorem for an asymptotically
relatively nonexpansive map defined on a proximal pair.
\section{Preliminaries}
\label{sec:1}
In this section, we introduce the notions of asymptotically relatively
nonexpansive maps and proximal parallel pairs having rectangle property. 
 Also, we give definitions and results related to this work.

\begin{definition} \cite{GK_Book}\label{P3:D:0}
A Banach space $X$ is said to have uniformly Kadec-Klee (UKK) norm
if and only if for any $\epsilon > 0$, there exists $\delta > 0$ such that
\begin{center}
   $\{x_n\} \subseteq B[0, 1]$, $x_n$ converges weakly to $x_0$,
  and ${\it sep}\{x_n\}:=\inf\{\|x_n - x_m\|: n \neq m\} > \epsilon$,
\end{center}
imply that
\begin{center}
  $\|x_0\| \leq 1 - \delta$.
\end{center}
\end{definition}

\begin{definition}\cite{GK_Book}
 A Banach space $X$ is said to be nearly uniformly convex (NUC) if and only if
  $X$ is a reflexive space having UKK norm.
\end{definition}

\begin{definition}\cite[4]{Eldr_2005}\label{P3:D:1}
Let $A$ and $B$ be nonempty
subsets of a Banach space $X$. The pair $(A, B)$ is said
to be a proximal pair if and only if for each $(x, y)$ in $A \times B$,
there exists $(x_1, y_1)$ in $A \times B$ such that
$\|x-y_1\| = {\it dist}(A, B) =\|y-x_1\|.$

In addition, if for each
$(x, y) \in A \times B,$ $(x_1, y_1) \in A \times B$ is a unique pair
  such that $\|x-y_1\| = {\it dist}(A, B)$ and $\|y-x_1\| = {\it dist}(A, B),$
then we say $(A, B)$ is a sharp proximal pair.
\end{definition}

\begin{definition}\cite{Espi_2008}\label{P3:D:2}
Let $A$ and $B$ be nonempty
subsets of a Banach space $X$.
The pair $(A, B)$ is said to
be a proximal parallel pair if and only if
\begin{enumerate}
  \item $(A, B)$ is a sharp proximal pair and
  \item there exists a unique $h \in X$ such that $B=A+h.$
\end{enumerate}
In this case, $a + h$ is called the proximal point in $B$ to $a$, for
any $a \in A$. 
\end{definition}

\begin{remark}\label{P3:R:0}
Let $(A, B)$ be a nonempty proximal pair in a
Banach space $X$. It was shown in \cite{Espi_2008} that
if $X$ is a strictly convex Banach space, then $(A, B)$ is a
proximal parallel pair.

In fact, the strict convexity of a Banach space is characterized as follows:

    $X$ is a strictly convex Banach space
    if and only if
   $(A, B)$ is a nonempty bounded closed convex
   proximal pair implies that there exists a unique $h \in X$ such that $B = A + h$.
\end{remark}

\noindent{\bf Proof} Suppose for every nonempty bounded closed convex
proximal pair $(A, B)$ in a
Banach space $X$ there exists a unique $h \in X$ such that $B = A + h$.
Now, it is claimed that $X$ is a strictly convex Banach space.

Suppose $X$ is not a strictly convex
Banach space. Then there exist distict points $x$ and $y$ in $X$
such that $\|x\| = 1 = \|y\|$ and $\|\frac{x+y}{2}\| = 1$.
It is also apparent that
$\|\alpha x+(1-\alpha)y\| = 1$, for every $\alpha \in ]0, 1[$.

Let $A = \{0\}$ and $B = \{\alpha x+(1-\alpha)y: \alpha \in [0, 1]\}$.
Then $(A, B)$ is a nonempty bounded closed convex
proximal pair in $X$. But $B \neq A + h$, for any
$h \in X$. Thus $X$ is a strictly convex space.\hfill{$\Box$}


\begin{definition}\cite{Goeb_1972}\label{P3:D:5}
Let $K$ be a nonempty subset of a Banach space $X$.
A mapping $T: K \rightarrow X$ is said to be asymptotically nonexpansive
if and only if there is a sequence $\{k_n\}$ in $[1, \infty)$ such that 
$\lim k_n = 1$ and 
\begin{center}
  $\|T^nx - T^ny\| \leq k_n\|x - y\|$,
\end{center}
for all $x, y \in K$. 
\end{definition}


\begin{definition}\label{P3:D:6}
Let $A$ and $B$ be nonempty subsets of a Banach space $X$.
A mapping $T$ from $A\cup B$ into $X$ is said to be asymptotically relatively
nonexpansive if and only if there is a sequence $\{k_n\}$ in $[1, \infty)$ such that 
$\lim k_n = 1$ and  
\begin{center}
  $\|T^nx - T^ny\| \leq k_n\|x - y\|$, 
\end{center}
for $x \in A$ and $y\in B$. 
\end{definition}

\begin{definition}\label{P3:D:7}
Let $(A, B)$ be a nonempty proximal parallel pair in a Banach space. Suppose
$T: A\cup B \rightarrow A\cup B$ is a map satisfying
$T(A) \subseteq B$ and $T(B) \subseteq A$. Define, for
$n \in \mathbb{N}$ and $x \in A\cup B$,
$u_n(x) = $ $\left\{
\begin{array}{ll}
        T^{n}(x'),  & \mbox{if}~~ n = 2m-1,\\
        T^{n}(x),  & \mbox{if}~~ n = 2m,
    \end{array}
 \right.$\\
 where $x' \in A\cup B$ such that $\|x - x'\| = \emph{dist}(A, B)$.
\end{definition}


\begin{lemma}\cite{Eldr_2006}\label{P3:L:0}
Let $A$ be a nonempty closed convex subset and $B$ be
a nonempty closed subset of a uniformly convex Banach space.
Let $\{x_n\}$ and $\{z_n\}$ be sequences in $A$ and $\{y_n\}$
be a sequence in $B$ satisfying:
\begin{enumerate}
  \item $\|x_n - y_n\|$ converges to ${\it dist}(A, B)$ and
  \item $\|z_n - y_n\|$ converges to ${\it dist}(A, B)$.
\end{enumerate}
Then $\|x_n - z_n\|$ converges to zero.
\end{lemma}

By observing Lemma \ref{P3:L:0}, the following concept was
introduced in \cite{Suzu_2009}.

\begin{definition}\cite{Suzu_2009}\label{P3:D:3}
A pair $(A, B)$ of nonempty subsets in a Banach space is
said to satisfy the property $UC$ if and only if the following holds:

If $\{x_n\}$ and $\{z_n\}$ are sequences in $A$ and $\{y_n\}$
is a sequence in $B$ such that
\begin{center}
  $\lim d(x_n, y_n) = {\it dist}(A, B)$
and $\lim d(z_n, y_n) = {\it dist}(A, B)$,
\end{center} then
$\lim d(x_n, z_n) = 0$.
\end{definition}

We introduce the following notion.

\begin{definition}\label{P3:D:4}
Let $(A, B)$ be a nonempty convex proximal parallel pair in
a Banach space $X$. The pair $(A, B)$ is said to have the
rectangle property if and only if
\begin{center}
  $\|x + h - y\| = \|y + h - x\|$, for any $x, y \in A$, 
\end{center}
where $h \in X$ such that $B = A + h$.

It is apparent that the proximal pair $(A, A)$ has the rectangle property.
\end{definition}

 \begin{remark}\label{P3:R:1}
Let $(A, B)$ be a nonempty proximal parallel pair in a
Banach space $X$.
It was proved in \cite{Kosuru_thesis} that if $X$ is a real Hilbert space,
then, for every $x, y \in A$ or $B$,
$x - y$ is orthogonal to $h$. 
\end{remark}

\noindent{\bf Proof} Let $x \in X$. It is known from \cite{Kham_Book} that if $K$ is a
nonempty closed convex subset of $X$, then
\begin{eqnarray}\label{P3:Eq:0}
  \langle z - P_K(x),~ P_K(x) -x\rangle &\geq& 0, ~for~ every~z~\in~ K,
\end{eqnarray}
where $P_K(x)$ is the best approximant of $x$ from $K$.

Now, let $x \in A$. Then $P_B(x) = x + h$ and hence
 it follows from eqn.(\ref{P3:Eq:0}) that
\begin{eqnarray}
  \langle (y+h) - (x+h),~ x+h -x\rangle &\geq& 0,~for~ every~ y+h ~\in~ B.
\end{eqnarray}\label{P3:Eq:0.1}
Similarly, let $x+h \in B$. Then $P_A(x+h) = x$. Thus
from eqn.(\ref{P3:Eq:0})
\begin{eqnarray}\label{P3:Eq:0.2}
  \langle y - x,~  x - (x+h)\rangle&\geq& 0,
\end{eqnarray}
for every $y \in A$. Equations (\ref{P3:Eq:0.1}) and (\ref{P3:Eq:0.2}) imply that
$\langle y - x, h\rangle = 0$, for $x, y \in A$. \hfill{$\Box$}

\begin{example}\label{P3:Ex:1}
Let $(A, B)$ be a nonempty proximal parallel pair in a real
Hilbert space. Then $(A, B)$ has the rectangle property.
\end{example}

\noindent{\bf Proof} It follows from Remark $\ref{P3:R:1}$ that
 $x - y$ is orthogonal to $h$, for every $x, y \in A$. Hence
$\|x - (y+h)\|^2 = \|x - y\|^2 + \|h\|^2 = \|x + h - y\|^2$.
Thus $(A, B)$ has the rectangle property. \hfill{$\Box$}

Next, we give an example of a proximal pair having the rectangle property
and the property $UC$ in a $NUC$ space.

\begin{example}\label{P3:Ex:2}
Consider the Day space
\begin{center}
$D_1 = \{x = (x^n): x^n \in \textit{l}_n ^1,
\|x\|_{D_{_1}}:=(\Sigma_{n \in \mathbb{N}}(\|x^n\|_1)^2)^{\frac{1}{2}} < \infty\}$,
\end{center}
where $\textit{l}_n ^1$ is the Banach space $\mathbb{R}^n$ with $\|.\|_{1}$.
The standard basis of $\mathbb{R}^n$ is denoted by $e_i ^n$,
for $1 \leq i \leq n$.
It is known from \cite{Bena_Book} that $D_1$ is a nearly uniformly convex
$(NUC)$ Banach space. Note that $D_1$ is not even strictly convex.

Consider, for $n \in \mathbb{N}$, $a_n = (a_n(i))_{i \in \mathbb{N}} \in D_1$, where
$a_n(i) = $ $\left\{
\begin{array}{ll}
        e_1 ^n \in \mathbb{R}^n, & \mbox{if}~~ i = n,\\
        0 \in \mathbb{R}^i,  & \mbox{if}~~ i \neq n.
    \end{array}
 \right.$\\
 Let $A = \overline{co}\{a_n : n \geq 2\}$ and $B = A + a_1$.
Suppose $x \in A$, then $x = (x(n)e_1 ^n)$, where $x(n)$ 
is a non-negative real number and $x(1) = 0$.
Similarly, if $y \in B$, then $y = (y(n)e_1 ^n)$, where
$y(n)$ is a non-negative real number and $y(1) = 1$. Now, for $x, y \in A$,
\begin{eqnarray*}
  \|x - y\|_{D_1} &=&(\Sigma_{i=2} ^{\infty}(\|(x(i) - y(i))e_1 ^i\|_1)^2)^{\frac{1}{2}}=(\Sigma_{i=2} ^{\infty}\mid x(i) - y(i)\mid^2)^{\frac{1}{2}}.
\end{eqnarray*}
Let $x , y \in A$. Then $x+ a_1$ and $y+ a_1$ are in $B$ and
\begin{eqnarray*}
\|x+a_1-y\|_{D_1} &=& (1 + \Sigma_{i=2} ^{\infty}(\|(x(i) - y(i))e_1 ^i\|_1)^2)^{\frac{1}{2}}.
\end{eqnarray*}
Hence
\begin{eqnarray}\label{P3:Eq:11}
  \|x+a_1-y\|_{D_1}&=& (1 + \|x - y\|_{D_1} ^2)^{\frac{1}{2}} = \|y + a_1 - x\|_{D_1}.
\end{eqnarray}
It follows from equation (\ref{P3:Eq:11}) that
$(A, B)$ is a proximal parallel pair having the rectangle property and
 ${\it dist}(A, B) = 1$.
 It is also easy to see from equation (\ref{P3:Eq:11}) that $(A, B)$ has the property $UC$.
\end{example}


\section{Main Results}
\label{sec:2}

 \indent We use the notation $u_n(.)$ as in Definition \ref
{P3:D:7}.



\begin{proposition}\label{P3:P:2}
Let $(A, B)$ be a nonempty proximal parallel pair having the 
property $UC$ in a Banach space $X$. Suppose
 $T: A\cup B \rightarrow A\cup B$ is an asymptotically relatively
nonexpansive map satisfying $T(A) \subseteq B$ and $T(B) \subseteq A$.
Then, for $x \in A$ and some subsequence $\{n_i\}$,
\begin{center}
  $u_{n_i}(x)$ converges weakly to $y$ if and only if
  $u_{n_i}(x+h)$ converges weakly to $y+h$,
\end{center}
where $y \in X$.
\end{proposition}

\noindent{\bf Proof} It is easy to see that
$\lim \|u_n(x) - u_n(x+h)\| = \emph{dist}(A, B)$, for all $x \in A$.
Since $(A, B)$ has the property $UC$,
$\lim\|u_n(x) + h - u_n(x+h)\| = 0$.

Suppose $u_{n_i}(x)$ converges weakly to $y$ for some $y \in X$.
Then the weak convergence of $u_{n_i}(x+h)$ to $y+h$ 
follows from
$\lim\|u_n(x)+h - u_n(x+h)\| = 0$.\hfill{$\Box$}

\begin{proposition}\label{P3:P:3}
Let $(A, B)$ be a nonempty bounded convex
proximal parallel pair in a Banach space $X$. Suppose
 $T: A\cup B \rightarrow A\cup B$ is an asymptotically relatively
nonexpansive map satisfying $T(A) \subseteq B$ and $T(B) \subseteq A$.
Further, assume that $(A, B)$ has the rectangle property and the property $UC$.
 Define $r_x: B \rightarrow \mathbb{R}$ and
$r_{y+h}: A \rightarrow \mathbb{R}$, respectively, by
$r_x(a+h) = \limsup_n \|u_{n}x - (a+h)\|$ and
$r_{y+h}(b) = \limsup_n \|u_{n}(y+h) - b\|$, for $x$ and $y \in A$.
Then
\begin{enumerate}
  \item the functions $r_x(.)$ and $r_{y+h}(.)$ are continuous and convex,
  \item $r_x(a + h) = r_{x+h}(a)$ and  hence
  $\inf_{y+h \in B}r_x(y+h) = \inf_{z \in A}r_{x+h}(z)$ and
  \item $r_x(u_n(a+h)) \leq k_nr_x(a+h)$, for every $n \in \mathbb{N}$.
\end{enumerate}
\end{proposition}

\noindent{\bf Proof}
1) It is easy to see that the functions
$r_x(.)$ and $r_{y+h}(.)$ are continuous and convex.

2) Let $x, a \in A$. Now, the rectangle property
of $(A, B)$ implies that
\begin{eqnarray*}
  \|a + h - u_n(x)\| &=& \|a - h - u_n(x)\|.\\
Thus,~\|a + h - u_n(x)\|&\leq& \|a - u_n(x+h)\| + \|u_n(x+h) - (h + u_n(x))\|.
\end{eqnarray*}
Since the pair $(A, B)$ has the property $UC$, $\limsup_n \|u_n(x+h) - (h + u_n(x))\| = 0$.
  Hence $r_x(a+h) \leq r_{x+h}(a)$. The other inequality, $r_{x+h}(a) \leq r_{x}(a+h)$,
follows exactly the same way. Thus 
$\inf_{y+h \in B}r_x(y+h) = \inf_{z \in A}r_{x+h}(z)$.

3) Let $n$ be an odd integer. Now, for $m \geq n$ and $m$ is
an odd integer
\begin{eqnarray*}
  \|u_m(x) - u_n(a+h)\|&=&\|T^m(x+h) - T^n(a)\|\\
  &\leq& k_n\|a - T^{m-n}(x+h)\|= k_n \|a - u_{m-n}(x+h)\|. 
  \end{eqnarray*}
In case that $m$ is an even integer,
\begin{eqnarray*}
  \|u_m(x) - u_n(a+h)\|&=&\|T^m(x) - T^n(a)\|\\
  &\leq& k_n\|a - T^{m-n}(x)\|= k_n \|a - u_{m-n}(x+h)\|. 
  \end{eqnarray*}
Since $r_{x+h}(a) = r_x(a+h)$,
$r_x(u_n(a+h)) \leq k_nr_x(a+h)$, for every $n = 2k-1$,
$k \in \mathbb{N}$. The other case, $n$ is an even integer, can be
proved similarly.\hfill{$\Box$}



%

\begin{lemma}\label{P3:L:1}
Let $(A, B)$ be a nonempty bounded closed convex proximal parallel pair in
a nearly uniformly convex $(NUC)$ Banach space $X$ and 
$T: A\cup B \rightarrow A\cup B$ be an asymptotically relatively
nonexpansive map satisfying $T(A) \subseteq B$ and $T(B) \subseteq A$.
Further, assume that $(A, B)$ has the rectangle property and the property $UC$.
Suppose $(K_1, K_2)$ is a subset of $(A, B)$ which is minimal with respect
to being nonempty, closed, and convex proximal pair satisfying:
\begin{eqnarray}\label{P3:Eq:0.5}
  x \in K_1 \Longrightarrow \omega_{w}(x) \subseteq K_1,
\end{eqnarray}
where $\omega_{w}(x)$ is the set of all limit points of the sequence $\{u_n(x)\}$
w. r. to the weak topology on $X$.
Then $r_x(y+h) = r_{y+h}(x)$,
for every $(x, y+h) \in K_1\times K_2$, where $r_x(.)$ and $r_{y+h}(.)$ are
as in Proposition \ref{P3:P:3}.

In this case, the constant $r_x(y+h) = r_{y+h}(x)$ is denoted by $r$.
\end{lemma}

\noindent{\bf Proof}
It is easy to see from Propostion \ref{P3:P:3} that the
functions $r_x(.)$, and 
$r_{y+h}(.)$ 
are lower semi continuous w. r. to the weak topology on $X$. Since $(K_1, K_2)$ is
a weakly compact convex pair, there exists $(x_0, y_0 + h)$ such that
$r_x: = r_x(y_0 + h)$  and
$r_{y+h} := r_{y+h}(x_0)$, where $r_x = \inf_{y+h \in K_2}r_x(y+h)$ and
$r_{y+h}= \inf_{x \in K_1} r_{y+h}(x)$.

Define $F_1:=\{y \in K_1 : r_{x+h}(y) \leq r_{x+h}\}$
and $F_2:=\{z \in K_2: r_x(z) \leq r_x\}$ for some $x \in K_1$. 
As the functions $r_x(.)$ and $r_{y+h}(.)$
are lower semicontinuous, the pair $(F_1, F_2)$ is nonempty.
Also, it follows from Proposition \ref{P3:P:3} that
$(F_1, F_2)$ is a proximal pair.

Now, it is claimed that $y \in F_1$ implies that $\omega_w(y) \subseteq F_1$.
Let $z \in \omega_w(y)$. Then
\begin{eqnarray*}
  r_{x+h}(z) &\leq& \liminf_{n_i} r_{x+h}(u_{n_i}y) \leq \limsup_n r_{x+h}(u_n(y))\\
   &\leq& \limsup_n k_n (r_{x+h}(y)) = r_{x+h}(y) = r_{x+h}.
 \end{eqnarray*}
Since $(K_1, K_2)$ is a minimal pair,
$r_x(y+h) = r_x = r_{x+h}(y)$, for every $y \in K_1$.

Let $z + h \in \omega_w(y+h)$. Then
\begin{eqnarray*}
  r_x=r_x(z+h) &=& \limsup_n \|z+h - u_n(x)\|\\
   &\leq& \limsup_n k_n \{\limsup_m\|x - u_{m-n}(y+h)\|\}\\
  &=&r_{y+h}(x) = r_{y+h}.
\end{eqnarray*}
The other inequality, $r_{y+h} \leq r_x$, can be
proved in a similar manner.

Hence $r_x(y+h) = r_{y+h}(x)$, for every $(x, y+h) \in K_1\times K_2$. \hfill{$\Box$}

\begin{lemma}\label{P3:L:3}
Let $(A, B)$ be a nonempty bounded closed convex proximal parallel pair in
a nearly uniformly convex $(NUC)$ Banach space $X$ and 
 $T: A\cup B \rightarrow A\cup B$ be a
continuous and asymptotically relatively
nonexpansive map satisfying $T(A) \subseteq B$ and $T(B) \subseteq A$.
Further, assume that $(A, B)$ has the rectangle property and the property $UC$.
Suppose $(K_1, K_2)$ is a subset of $(A, B)$ which is minimal with respect
to being nonempty, closed, and convex proximal pair satisfying:
\begin{enumerate}
  \item $x \in K_1 \Longrightarrow \omega_{w}(x) \subseteq K_1$, where
  $\omega_w(x)$ is as in Lemma \ref{P3:L:1},
  \item for each $x \in K_1$, every subsequence of $\{u_n(x)\}$ admits a
  norm convergent subsequence and
  \item for each $x \in K_1$, $\omega(x)$ is a norm compact set, where
  $\omega(x)$ is the set of all limit points of the sequence $\{u_n(x)\}$ w. r. to
  the norm topology.
  \end{enumerate}
Then there exists $x \in K_1$ such that
$\|x - Tx\| = \emph{dist}(A, B)$.
\end{lemma}

\noindent{\bf Proof} It is to be noted that the property $UC$ of $(A, B)$
and Proposition \ref{P3:P:2} imply
that $\omega_w(x+h) = \omega_w(x) + h$ and
$\omega(x+h) = \omega(x) + h$, for every $x \in A$.

If $K_1 = \{x_0\}$, then $u_n(x_0)$ converges to $x_0$ and
hence the continuity of $T$ implies
that $u_n(x_0 + h)$ converges to $T(x_0)$.
Since $\omega(x_0 + h) = \omega(x_0) + h$, $T(x_0) = x_0 + h$.

It is claimed that $\delta(K_1) = 0$. Suppose $\delta(K_1) > 0$.

Fix $x \in K_1$, define $F_1 := \omega(x)$ and $F_2 := F_1 + h$.
Then $(F_1, F_2)$ is a proximal parallel pair. As $T$ is continuous,
$(T(F_2), T(F_1)) \subseteq (F_1, F_2)$. Let
$\mathcal{F}$ be the set of all nonempty and closed proximal pair
subsets $(H_1, H_2)$ of $(F_1, F_2)$ which satisfies:
\begin{center}
  $x \in H_1 \Longrightarrow \omega(x) \subseteq H_1$.
\end{center}
Define a relation on $\mathcal{F}$ by
$(H_1, H_2) \leq (D_1, D_2)$ if and only if
$(H_1, H_2) \subseteq (D_1, D_2)$.

It is easy to see that every nonempty totally ordered
subset $\mathcal{T}$ of $\mathcal{F}$ has a lower bound.
Hence by Zorn's lemma, we obtain a minimal
proximal pair $(H_1, H_2)$ which satisfies
$(T(H_2), T(H_1)) = (H_1, H_2)$, that is,
$(u_1(H_1), u_1(H_2)) = (H_1, H_2)$ and hence, for every $n \in \mathbb{N}$,
$(u_n(H_1), u_n(H_2)) = (H_1, H_2)$.

Now, it is easy to see that $H_1$ is a singleton set.
For, since $(H_1, H_2)$ is a compact proximal parallel pair,
$(\overline{co}H_1, \overline{co}H_2)$ is a compact
convex proximal parallel pair and hence there exists
$(x, y+h) \in (\overline{co}H_1, \overline{co}H_2)$ such that
\begin{center}
$\delta(x, H_2) < \delta(H_1, H_2)$ and
$\delta(y+h, H_1) < \delta(H_1, H_2)$.
\end{center}
Define $D_1 = \{y \in H_1 : \delta(y, H_2) \leq \alpha\}$ and
$D_2 = \{z \in H_2 : \delta(z, H_1) \leq \alpha\}$, where $\alpha = \delta(x, H_2)$.
Since $(A, B)$ has the rectangle property,
$(D_1, D_2)$ is a nonempty closed 
 proximal parallel pair. Also, it is apparent that
$(D_1, D_2) \subsetneqq (H_1, H_2)$.

Now, it is claimed that
$x \in D_1 \Longrightarrow \omega(x) \subseteq D_1$.
Let $x \in D_1$ and $y = \|.\|-\lim u_{n_i}(x)$ for some
subsequence $\{n_i\}$.

Since $(u_n(H_1), u_n(H_2)) = (H_1, H_2)$ for every $n \in \mathbb{N}$,
 there exists $x_n \in H_1\cup H_2$ such that $u_n(x_n) = x$ for every $x \in H_1\cup H_2$.
Now, let $z \in H_2$
\begin{eqnarray*}
  \|y - z\| & = & \lim \|u_{n_i}(x) - z\| \leq \limsup \|u_n(x) - u_n(z_n)\|\\
  &\leq& \limsup k_n\|x - z_n\| \leq \delta(x, H_2).
\end{eqnarray*}
Hence $\delta(y, H_2) \leq \delta(x, H_2)$ and $y \in D_1$.
This shows that $\delta(H_1) = 0$ and
so is $K_1$. \hfill{$\Box$}

\begin{theorem}\label{P3:T:1}
  Let $(A, B)$ be a nonempty bounded closed convex proximal parallel pair in
a nearly uniformly convex $(NUC)$ Banach space. Suppose
$T: A\cup B \rightarrow A\cup B$ is a
continuous and asymptotically relatively
nonexpansive map satisfying $T(A) \subseteq B$ and $T(B) \subseteq A$.
Further, assume that $(A, B)$ has the rectangle property and the property $UC$.
Then there exists $x \in A\cup B$ such that
$\|x - Tx\| = {\it dist}(A, B)$.
\end{theorem}

\noindent{\bf Proof}
In view of Lemma \ref{P3:L:3}, it is enough to prove that
$(A, B)$ has a minimal proximal parallel pair
$(K_1, K_2)$ which satisfies:
\begin{enumerate}
  \item $x \in K_1 \Longrightarrow \omega_{w}(x) \subseteq K_1$, where
  $\omega_{w}(x)$ is the set of all limit points of the sequence $\{u_n(x)\}$
  w. r. to the weak topology on $X$,
  \item for each $x \in K_1$, every subsequence of $\{u_n(x)\}$ admits a
  convergent subsequence w. r. to the norm topology and
  \item for each $x \in K_1$,
  $\omega(x)$ is a norm compact set.
\end{enumerate}
Let $\mathcal{F}$ be the set of all nonemtpy closed convex
proximal parallel pair subset $(K_1, K_2)$ of $(A, B)$ satisfying:
\begin{center}
  $x \in K_1 \Longrightarrow \omega_w(x) \subseteq K_1$
\end{center}
Define a relation $\leq$ on $\mathcal{F}$ by
$(H_1, H_2) \leq (K_1, K_2)$ if and only if $(H_1, H_2) \subseteq (K_1, K_2)$.
It is easy to see that every nonempty totally ordered subset of
$\mathcal{F}$ has a lower bound. Hence by Zorn's Lemma
$\mathcal{F}$ has a minimal closed convex proximal parallel pair $(K_1, K_2)$
satisfying $x \in K_1 \Longrightarrow \omega_w(x) \subseteq K_1$.

It is claimed that $K_1 = \{x\}$ for some $x \in A$.

If $r = d$, where $r = r_x(y+h) = r_{y+h}(z)$
($r_x(.)$ and $r_{y+h}(.)$ are as in Proposition \ref{P3:P:2}), for 
$x, y$ and $z$ in $K_1$ and
$d:= {\it dist}(A, B)$, then, for every $x \in K_1$,
$u_n(x)$ converges to $y$ for every $y \in K_1$.
Hence $K_1 = \{x\}$ for some $x \in A$ and
the continuity of $T$ implies that $Tx = x + h$. 

Now, suppose that $r > d$ 
 and, for some $x \in K_1$,
there exists a sequence $\{m_i\}$ such
that $u_{m_i}(x)$ does not have any norm convergent subsequence.
As $\{u_{m_i}(x)\}$ is not compact, there exists $\epsilon > 0$
and a sequence $\{l_i\}$ such that
$\|u_{m_{_{l_i}}}(x) - u_{m_{_{l_j}}}(x)\| \geq \epsilon$, for every $i \neq j$.

Let $y_{_i} = u_{m_{_{l_i}}}(x)$. Since $K_1$ is weakly compact,
$\{y_{_i}\}$ has a weakly convergent subsequence. 
We assume that $y_{_i}$ converges weakly to $y_{_0}$ for some
$y_{_0} \in K_1$.

Since $X$ is a NUC space, there exists $\delta \in~ ]0, 1[$ corresponding
to $\frac{\epsilon}{2r}$ such that $r(1 - \delta) < r$. Now,
choose $\epsilon_1\in~ ]0, r[$ such that
$r_1: = (r + \epsilon_1)(1 - \delta) < r.$

Since $r_{x}(x+h) = r$ and $\lim k_n = 1$,
there exists $N \in \mathbb{N}$ such that
\begin{center}
 $\|x+h - u_n(x)\| \leq r + \frac{\epsilon_1}{2}$ and
 $(k_n - 1)(r + \frac{\epsilon_1}{2}) \leq \frac{\epsilon_1}{2}$,
 for all $n \geq N$.
\end{center}
Let $n \geq N$ be fixed. Now, for every $i \geq n + N$,
\begin{eqnarray*}
  \|y_{_i} - u_n(x+h)\|&=&\|u_{m_{_{l_i}}}(x) - u_n(x+h)\|\\
  &\leq& k_n\|x+h - u_{m_{_{l_i}}-n}(x)\| \leq r+\epsilon_1.
\end{eqnarray*}
Let $z_{_i} = \frac{y_{_{i}} - u_{n}(x+h)}{r+\epsilon_1}$,
for every $i \geq n + N$. Then $z_{_i} \in B[0, 1]$,
$z_{_i}$ converges weakly to $\frac{y_{_0} - u_{n}(x+h)}{r+\epsilon_1}$,
for every fixed $n \geq N$ and ${\it sep}\{z_{_i}\} \geq \frac{\epsilon}{r+\epsilon_1}$.

Hence $\|y_{_0} - u_n(x+h)\| \leq (r+\epsilon_1)(1 - \delta)$,
for all $n \geq N$.
Therefore,
\begin{center}
  $r:=\limsup_n\|y_{_0} - u_{n}(x+h)\| \leq (r+\epsilon_1)(1 - \delta) < r$.
\end{center}
This contradiction shows that
every subsequence of $\{u_n(x)\}$ admits a
norm convergent subsequence, for every $x \in K_1$. 

Now, it is claimed that, for every $x \in K_1$,
$\omega(x)$ is a norm compact set.

Suppose that $\omega(x)$ is not
a norm compact set for some $x \in K_1$. Then there exists a sequence $\{y_{_n}\}$
in $\omega(x)$ such that ${\it sep}\{y_{_n}\} \geq \epsilon$ for some $\epsilon > 0$.

Since $y_{_k} \in \omega(x)$ for every $k \in \mathbb{N}$,
there exists $n_k \in \mathbb{N}$ such that
$\|y_{_k} - u_{n_{_k}}(x)\| \leq \frac{\epsilon}{3}$.
Note that $n_k$ can be chosen such that $n_k < n_{k+1}$ for every $k \in \mathbb{N}$.
Define $z_{_k} = u_{n_{_k}}(x)$. Then from triangle inequality it follows that
${\it sep}\{z_{_k}\} \geq \frac{\epsilon}{3}$. Thus 
$\{z_{_k}\}$ does not admit any norm convergent subsequence,
a contradiction to 
every subsequence of $\{u_n(x)\}$ admits a
norm convergent subsequence. \hfill{$\Box$}

The following example illustrates Theorem \ref{P3:T:1}.

\begin{example}
Consider the series $\Sigma_{n=1} ^{\infty} (-1)^{n+1}\frac{1}{n}$,
 define $b_n := \frac{1}{2n-1}-\frac{1}{2n} > 0$. Then
$\Sigma_{n=1} ^{\infty} b_n$ converges to $log_{e} 2$.
Now, let $c_{n+1} = (1/e^{b_n})$, for $n \in \mathbb{N}$. Then $0 < c_n < 1$ and
$\Pi_{n=2} ^{\infty}c_n = \frac{1}{2}$.

Let $a_n = (a_n(i))_{i \in \mathbb{N}} \in D_1$, where
$a_n(i) = $ $\left\{
\begin{array}{ll}
        e_1 ^n \in \mathbb{R}^n, & \mbox{if}~~ i = n,\\
        0 \in \mathbb{R}^i,  & \mbox{if}~~ i \neq n.
    \end{array}
 \right.$\\
 Let $A = \overline{co}\{a_n : n \geq 2\}$, and $B = A + a_1$.
 In Example \ref{P3:Ex:2}, it is shown that $(A, B)$
 is a bounded closed convex proximal parallel pair having the
 rectangle property and the property $UC$
 in the nearly uniformly convex $(NUC)$ Banach space $D_1$.

 Define $T: A\cup B \rightarrow A\cup B$ by
$T(x) = $ $\left\{
\begin{array}{ll}
        (1, 0, x(2)^2e_1 ^3, c_2(x(3)e_1 ^4), ...),  & \mbox{if}~~ x \in A,\\
        (0, 0, x(2)^2e_1 ^3, c_2(x(3)e_1 ^4), ...),  & \mbox{if}~~ x \in B.
\end{array}
 \right.$\\
It is clear that $T$ is a continuous map satisfying
$T(A) \subseteq B$ and $T(B) \subseteq A$. Also,
for $(x, y) \in A\times B$,
\begin{eqnarray*}
  \|Tx - Ty\|_{D_1}&=&(1+\|(x(2) ^2-y(2) ^2)e_1 ^3\|_1 ^2+\|c_2(x(3)-y(3))e_1 ^4\|_1 ^2 + ...)^{\frac{1}{2}}\\
  &\leq& 2(1 + \mid x(2)-y(2) \mid^2 + \mid x(3)-y(3)\mid^2 + ...) = 2\|x - y\|_{D_1}.
\end{eqnarray*}
Similarly, it can be shown that
$\|T^n(x) - T^n(y)\|_{D_1} \leq 2(\Pi_{i=2} ^n c_i) \|x - y\|_{D_1}$,
 for $n \geq 2$ and $(x, y) \in A\times B$.
Thus, $T$ is an asymptotically relatively nonexpansive map. Hence from
 Theorem \ref{P3:T:1} it follows that $T$ has a best proximity point in $A\cup B$.

Indeed, $(0, (0, 0), (0, 0, 0),...) \in A$ is a best proximity point of $T$. \hfill{$\Box$}
\end{example}

\begin{theorem}\label{P3:T:2}
  Let $(A, B)$ be a nonempty bounded closed convex proximal parallel pair in
a nearly uniformly convex $(NUC)$ Banach space. Suppose
 $T: A\cup B \rightarrow A\cup B$ is a
continuous and asymptotically relatively
nonexpansive map satisfying $T(A) \subseteq A$ and $T(B) \subseteq B$.
Further, assume that $(A, B)$ has the rectangle property and the property $UC$.
Then there exist $x \in A$ and $y\in B$ such that $Tx = x$, $Ty = y$ 
and $\|x - y\| = \emph{dist}(A, B)$.
\end{theorem}

\noindent{\bf Proof} This can be proved in a similar way
 as Theorem \ref{P3:T:1}. \hfill{$\Box$}

If $A = B$, then we obtain the next result, which is Corollary 1 in \cite{Xu_1991},
from Theorem \ref{P3:T:1} and
Theorem \ref{P3:T:2}.

\begin{corollary}\cite{Xu_1991}
  Let $X$ be a nearly uniformly convex Banach space and
  $A$ be a nonempty bounded closed convex subset of $X$.
  Suppose $T:A \rightarrow A$ is an asymptotically nonexpansive
  map. Then $T$ has a fixed point in $A$.
\end{corollary}

It is known that every uniformly convex Banach space is a nearly uniformly convex (NUC) 
space. Hence we obtain the following result.

\begin{corollary}\cite{Goeb_1972}
 Let $A$ be a nonempty bounded closed convex subset of
 a uniformly convex Banach space $X$ 
 and $T:A \rightarrow A$ be an asymptotically nonexpansive
 map. Then $T$ has a fixed point in $A$.
\end{corollary}

\noindent{\bf Acknowledgements} 
The authors thank the referees for the suggestions and comments. 
The first author, S. Rajesh,
thanks the University Grants Commission (India) for the financial support
provided as a form of Research Fellowship to carry out this research
work at Indian Institute of Technology Madras, Chennai.

\end{document}